\documentclass[11pt]{amsart}
\usepackage{latexsym}
\usepackage{amsthm}
\usepackage{amsmath}
\usepackage{amsthm}  
\usepackage{epsfig}
\usepackage{amsfonts}
\usepackage{amssymb}
\usepackage{epsfig}
\usepackage{graphicx}
\usepackage{graphics}

\newtheorem{theorem}{Theorem}

\newtheorem*{definition}{Definition} 
\newtheorem{proposition}[theorem]{Proposition}
\newtheorem{corollary}[theorem]{Corollary}
\newtheorem*{example}{Example} 
\newtheorem*{remark}{Remark} 

\def\R{\mathbb R}

\def\p{\mathbf{p}}

\def\x{\mathbf{x}}
\def\y{\mathbf{y}}

\def\0{\mathbf{0}}
\def\al{\mathbf{\alpha}}
\def\rhat{\hat\rho}
\if@amsmath

\else

\fi

\begin{document}

\title{M\"obius transformations of polygons and partitions of 3-space}
\author{Richard Randell}
\address{Department of Mathematics, University of Iowa,Iowa City IA 52242.}
\author{Jonathan Simon}
\author{Joshua Tokle}

\thanks{Research partially supported by NSF Grant DMS\,0107209. Email
richard-randell@uiowa.edu or jonathan-simon@uiowa.edu.}

\maketitle

\today

\begin{abstract}

The image of a polygonal knot $K$ under a spherical inversion of
$\mathbb{R}^{3} \cup\infty$ is a simple closed curve made of arcs of
circles, perhaps some line segments, having the same knot type as
the mirror image of $K$. But suppose we reconnect the vertices of
the inverted polygon with straight lines, making a new polygon
$\widehat{K}$. This may be a different knot type. For example, a
certain 7-segment figure-eight knot can be transformed to a
figure-eight knot, a trefoil, or an unknot, by selecting different
inverting spheres. \emph{Which knot types can be obtained from a
given original polygon $K$ under this process?} We show that for
large $n$, most n-segment knot types cannot be reached from one
initial n-segment polygon, using a single inversion or even the
whole M\"{o}bius group.

The number of knot types is bounded by the number of complementary
domains of a certain system of round 2-spheres in $\mathbb{R}^{3}$.
We show the number of domains is at most polynomial in the number of
spheres, and the number of spheres is itself a polynomial function
of the number of edges of the original polygon.  In the analysis, we obtain an exact formula for the number of complementary domains of any collection of round $2$-spheres in $\R^3$.  On the other hand,
the number of knot types that can be represented by $n$-segment
polygons is exponential in $n$.

Our construction can be interpreted as a particular instance of
building polygonal knots in non-Euclidean metrics. In particular,
start with a list of $n$ vertices in $\mathbb{R}^{3}$ and connect
them with arcs of circles instead of line segments: Which knots can
be obtained? Our polygonal inversion construction is equivalent to
picking one fixed point $p \in \R^3$ and replacing each edge of $K$
by an arc of the circle determined by $p$ and the endpoints of the
edge.

\end{abstract}

\section{Introduction}

Inversion of 3--space through a sphere is a well-known transformation of
$\mathbb{R}^{3} \cup\mathbf{\infty}$. If $S_{\mathbf{p}, r}$ is the round
sphere of radius $r$ centered at the point $\mathbf{p}$, the mapping
\[
\rho(\mathbf{x}) =\mathbf{p} + \frac{r^{2}(\mathbf{x}-\mathbf{p})}%
{|\mathbf{x}-\mathbf{p}|^{2}}%
\]
sends $\mathbf{p} \leftrightarrows\mathbf{\infty}$ and fixes the points of
$S_{\mathbf{p} , r}$. Sphere inversions are conformal maps, take circles and
lines to circles or lines, spheres and planes to spheres or planes. We
consider planes as spheres through infinity, and include reflections in planes
as inversions through spheres. The composition of an even number of sphere
inversions [is orientation-preserving, and] is called a M\"obius
transformation. The inversions form a group under composition. For the
purposes of this paper, we want to include orientation-reversing maps as well,
so we will use the term M\"obius transformation to mean any composition (even
or odd) of sphere inversions. (See e.g. Section 2.2 of \cite{ThurstonLevyBook}
for a general introduction to M\"obius transformations and their classical interpretations.)

\begin{figure}[h]

 \includegraphics[bb=0 0 230 230]{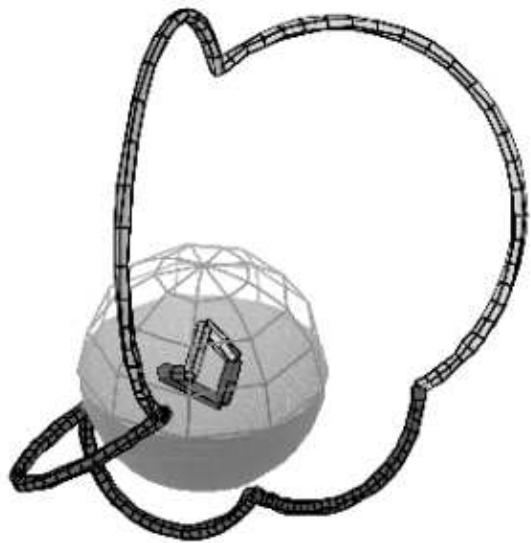}
\caption{Inverting a polygon through a sphere to get a
circle--arc polygon}
\label{Fig1SimpleInversion}
\end{figure}

We denote the mirror image of a knot $K$ by $K^{*}$. 
If $K$ is a polygonal knot in $\mathbb{R}^{3}$, and $\rho$ is an inversion
(whose center does not lie on $K$), the set $\rho(K)$ is a closed curve
(see Figure \ref{Fig1SimpleInversion}) that is the union of arcs of circles and
perhaps (if the center of inversion is colinear with some
edges of $K$) some line segments. The knot $\rho(K)$ is isotopic to $K^{*}$.

\begin{figure}[h]

\includegraphics[bb=0 0 178 178]
{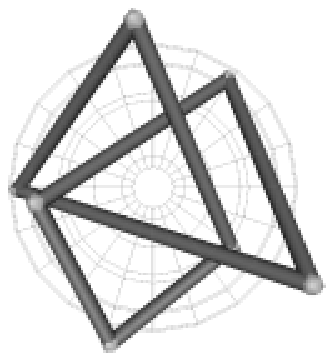}
\vspace{0.1in}

\includegraphics[bb=0 0 178 178]
{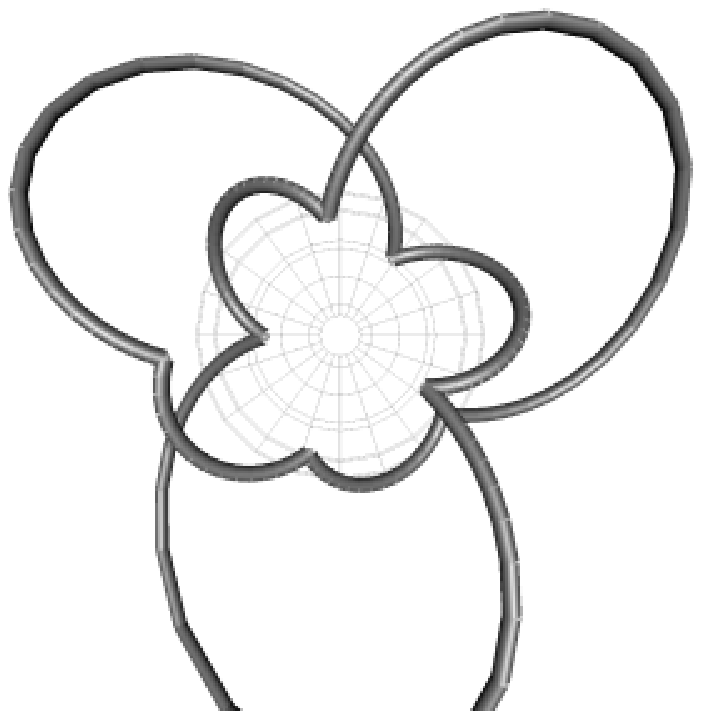}
\vspace{0.1in}
\includegraphics[bb=0 0 178 178]
{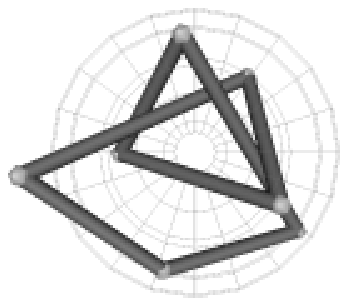}
\hfill
\caption{For this inversion of the upper  trefoil $K$, $\rho(K)$ is a  LH trefoil, and  $\rhat(K)$ is a RH trefoil.}
\label{TrefBothInversions}
\end{figure}

\begin{figure}[h]

\includegraphics[bb=0 0 178 178]
{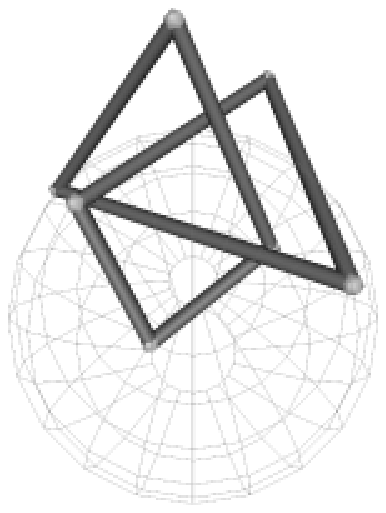}
\vspace{0.1in}

\includegraphics[bb=0 0 178 178]
{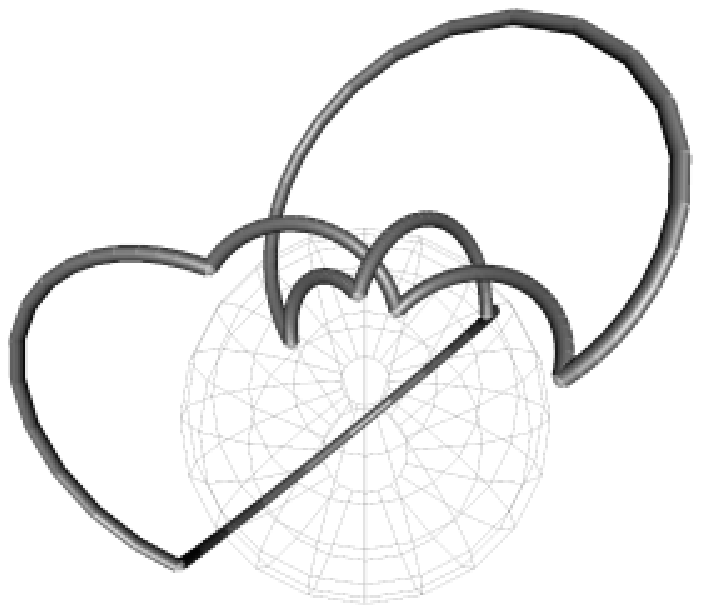}
\vspace{0.1in}
\includegraphics[bb=0 0 178 178]
{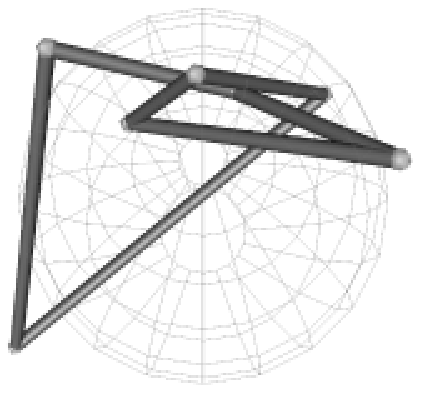}
\hfill
\caption{With the same polygon $K$ and a different inversion sphere, we have $\rho(K)$ is a  LH trefoil, and  $\rhat(K)$ is an unknot.}
\label{TrefDifferentSphere}
\end{figure}

\emph{What would happen if instead of using the spherically inverted edges of $K$ to
connect the vertices of $\rho(K)$, we formed a new polygon by connecting the
successive vertices of $\rho(K)$ with line segments?} Call this  operation {\em polygonal inversion} and denote the resulting polygon by
$\hat\rho(K)$.  (See Figures \ref{TrefBothInversions}, \ref{TrefDifferentSphere}).  The polygon $\hat\rho(K)$ may be singular, but in general it is an embedded knot, perhaps very different from $K^*$:  {\em Which knot types can be obtained this way from a given starting polygon $K$}?

\begin{figure}[h]
\centering

\includegraphics[bb=0 0 230 230]
{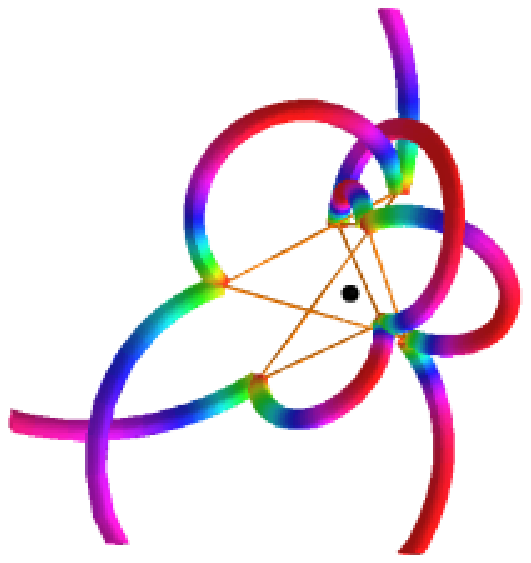}

\vspace{0.1in}

\includegraphics[bb=0 0 230 230]
{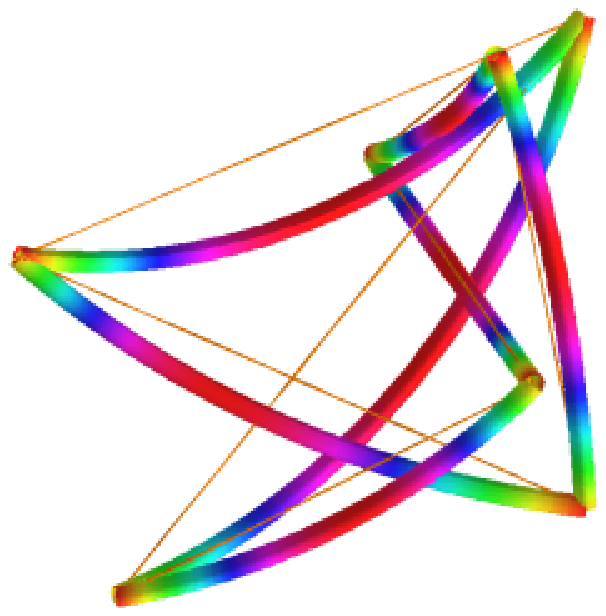}

\caption{ Inversions of the same figure-eight knot, giving LH and RH trefoils, here seen as arcs of circles through the center of inversion. (In the first figure, the center is near the knot and some of the circles are very large.  In the second figure,  the center is far away, so the circle-arcs are closer to the polygon edges.)}

\label{CircleArcPictures}
\end{figure}

Here is another way to visualize the knot type of $ \hat\rho(K)$.
Keep the vertices of $K$ the same and replace each edge by a certain circle-arc, as follows (Figure \ref{CircleArcPictures}):  For each edge of $K$, construct the circle
containing the endpoints $\x,\y$ of the edge and the inversion center $\p$.
Replace the edge $[\x \y]$ of $K$ by the arc $\langle \x \y \rangle$ of that circle not containing $\mathbf{p}$;
call the circular polygon we get $\widetilde{K}$.   See Figure \ref{CircleArcPictures}.
Then $\widetilde{K} =
\rho(\hat\rho(K))$, so $\widetilde{K}$ is the same knot type as $\hat\rho(K)^*$.
Seen this way, our question of which knot types can arise as $\widehat{K}$ is
a special case of the problem of determining which knot types can arise if one
replaces the edges of a polygonal knot by a some set of circle-arcs, passing in this way
from Euclidean polygons to non-Euclidean polygons.

We originally were led to the study of inverting polygons from questions about
knot energies, in particular to find polygons that look very different
geometrically yet have the same, or at least close, polygonal knot energies.
Suppose $K$ is an inscribed fine-mesh approximation of a smooth curve $C$.
According to \cite{RawdonSimon}, the \emph{minimum distance energy}
$E_{md}(K)$ is close to the \emph{M\"obius energy} $E_{O}$ of $C$. If we apply
a sphere inversion $\rho$, then $E_{O}(\rho(C)) = E_{O}(C)$, even though they
may look very different in Euclidean geometry. The polygon $\hat\rho(K)$ will
be an inscribed polygon for $\rho(C)$. If the center $p$ is not too close to
$K$, then $\hat\rho(K)$ will be a close approximation of $\rho(C)$ and so
$E_{md}(\hat\rho(K)) \approx E_{O}(\rho(C)) = E_{O}(C) \approx E_{md}(K)$.
This approach has led us to many local minima for $E_{md}$ with similar
energies.

As just noted, if $K$ looks nearly smooth, then $\hat\rho(K)$ may also look
nearly smooth and be of the same knot type as $K$ (in particular, isotopic to $K^{*}$) But
for knots made of relatively few edges, the knot type of $\hat\rho(K)$ may vary.

In {\bf Section \ref{Examples}}, we give several examples of polygonal inversions: a right-handed trefoil that can be changed to itself, or an unknot, or a LH trefoil; and  a 7-segment figure-eight knot where the polygons $\hat\rho(K)$
include an unknot, a right-handed trefoil knot, and a left-handed trefoil
knot, along with figure-eights. These are  (\cite{Randell}  \cite{Negami}  \cite{JinandKim})  all the knot types that can be represented as a polygon with 7 edges.
\emph{Are there other situations where a particular starting polygon is
``universal" like this particular 7-segment figure-eight?}

Note that if $K$ can be inverted to $\hat\rho(K)$, then the same polygonal inversion takes $\hat\rho(K)$ to $K$; so polygonal inversion can make knots more complicated as well as simpler.  Also if we allow compositions of inversions, then for any ``universal" polygon $K$, each of the polygonal inverses $\hat\rho(K)$ also is universal.

Finally, at the other extreme, we describe a large class of polygons (representing all knot types) for which polygonal inversions yield only the original knot type and its mirror image.
 \vspace{0.2in}
 
We show in this paper that, in general, the knot types that can be obtained
from a given original $n$-segment knot $K$ are a small portion of all the
$n$-segment knots.  
\begin{theorem}
There are exponentially many $n$-segment
knot types, but polygonal inversions acting on a given polygon $K$ can only yield
polynomially many knot types.
\end{theorem}


For explicit bounds, see Theorem \ref{maximal}, Corollary \ref{BoundKnots}, and Theorem  \ref{ExponentiallyMany}.  In particular, no 72-edge polygon can be inverted to reach all knot types that can be formed with 72 edges.  Since the number of $n$-edge knot types is much larger than the known minimum bound that we use in our calculations, we expect that the critical number of edges is considerably smaller than 72.

\vspace{0.1in}

In {\bf Section \ref{Basics}}, we give some basic properties of sphere inversions.
In {\bf Section \ref{ReduceToSpheres}} we show that the number of knot types arising
from a single inversion of a given $n-$edge polygon $K$ is bounded by the number of
complementary domains of a certain system of $m$ round 2-spheres (including
the possibility that some are planes), where $m \leq\ n(n-3)/2$.

In {\bf Section \ref{Domains}} we show that a system of $m$ round 2-spheres
(including perhaps planes) in $\R^3$ has at most $2\binom{m}{3}+2m$
complementary domains. This is a key part of our overall argument, and may be
of independent interest. We actually obtain an exact formula for the number of domains based on the intersection pattern of the spheres.  The proof is complicated by the fact that the spheres
may not be mutually transversal. Combined with the previous section, this
says that the number of knot types obtainable from $K$ by one inversion is at
most on the order of $n^{6}$.

In {\bf Section \ref{ExponentiallyMany}}, we show that the number of $n$-segment
knot types grows exponentially with $n$. Thus, for large $n$, the number of knot types that can be reached by inversion from a given
polygon $K$ is much smaller than the total number of n--segment knot types. In particular, for large $n$, there cannot exist an $n$-segment knot that is universal like the 7-segment figure-eight.

Finally, in {\bf Section \ref{MobiusGroup}}, we show that if we let the whole group of M\"obius transformations act on a given polygon, that is compose arbitrary numbers of sphere inversions before reconnecting the vertices, the number of knot types arising is no larger than twice the number from just single spherical inversions.  Thus our theorem
(use 75 edges instead of 72) applies as well to the
action of the whole M\"obius group.

\section{Examples of changing knot type via inversions}

\label{Examples}

\subsection{M\"obius inversion vs. polygonal inversion}
In Figure \ref{TrefBothInversions} and Figure \ref{TrefDifferentSphere} we show a particular 7-segment RH-trefoil knot and two different inversion spheres. The  circle-arc ``polygons" $\rho(K)$ must be equivalent to $K^*$.  However, in the first example, $\rhat(K)$ is a RH-trefoil, while in the second example, $\rhat(K)$ is unknotted. 
For sufficiently large numbers $c$, if we use a sphere of radius $c$ centered at the point $(c,0,0)$, then the polygon $\rhat(K)$ is be equivalent to $K^*$ (see Section \ref{Basics}).

\subsection{A 7-segment figure-eight knot that is ``universal"}

For a certain 7-segment figure-eight knot, there are spheres of inversion that
yield unknots,  left-handed trefoils, right-handed trefoils and figure-eight knots. 
These are  (\cite{Randell}  \cite{Negami}  \cite{JinandKim})  all the knot types that can be
represented as polygons with $\leq7$ edges.

The vertex matrix of the initial figure-eight $K$ is:

\begin{displaymath}
\left(
\begin{array}{cccccc}

  -1  & \,&  -13& \,&     24\\
   -9  & \,&   24 & \,&     19\\
  -27  & \,&  -15  & \,&  -20\\
   45   & \,&  3& \,&    -2\\
  -23  & \,&    7 & \,&    34\\
   30  & \,&  -15 & \,&   -37\\
  -16   & \,&   10 & \,&   -17\\

\end{array}
\right)
\end{displaymath}

By choosing different centers of inversion (the radii do not matter; see first paragraph of Section \ref{ReduceToSpheres}) we obtain
\begin{displaymath}
\begin{array}{ccc}
\text{Center of Inversion} &\,& \text{ Knot Type of } \hat\rho(K)\\
\cline {1-1} \cline{3-3}
(0,0,0) &\,& \text{Unknot} \\
(-6, -6, -6) &\,& \text{RH trefoil} \\
(100,100,100) &\,& \text{LH trefoil} \\
(1000,1000,1000) & \, & \text{mirror image } K^*

\end{array}
\end{displaymath}

In addition to studying how the operation $\hat\rho$ changes topological knot type, one might also ask how it changes {\em polygonal} knot type, that is equivalence in which polygons must be deformed keeping the number of edges fixed.  It is shown in \cite{Calvo} that any seven-segment figure-eight is equivalent, in this strong sense, to its mirror image.  From this, and our analysis in Section \ref{ReduceToSpheres}, it follows that for far-away centers of inversion, $\hat\rho(K)$ is polygonally equivalent to $K$, not to the reverse of $K$ (which is polygonally different, again by \cite{Calvo}).

\subsection{Leaving a knot fixed}

Suppose $K$ is constructed so that its vertices all lie on one round sphere $S$.  Every knot  type can be realized this way.  (For example, given a tame knot type, there exists a number $N$ such that \cite{Negami} {\em every} linear embedding in $\R^3$ of the complete graph on $N$ points contains a cycle of the given knot type).   Inversion in $S$ takes $K$ to a
circle-polygon $\rho(K)$ of type $K^*$; on the other hand, $\rhat(K)= K$ (all vertices are fixed by the inversion).
In fact (see Section \ref{ReduceToSpheres}), for such $K$, all spheres centered inside $S$ give $\rhat(K)\approx K$, while all spheres centered outside $S$ give $\rhat(K)\approx K^*$.

\subsection{General remarks}
 \label{ackn}
Note that $\rhat$ is an inversion in the sense that for any polygon $K$, $\rhat(\rhat(K))=K$.  Thus if some $\rhat$ takes an $8_{19}$ polygon to an unknot, then the same $\rhat$ takes the unknotted polygon to the $8_{19}$.  So polygonal inversion can make knots more complicated as well as simpler.

 Sometimes an inversion distorts a
polygon so much that it is hard to see the knot type just by looking at the
output polygon. To confirm experiments, after inverting the polygons, we have
relaxed them using the minimum-distance energy $U_{md}$ and software
\emph{MING} \cite{WuMing}. We also used a modification of MING done by E. Rawdon, which incorportates sphere-inversions.  More recently, we have used R. Scharein's {\em KnotPlot} and an implementation of MING ported into KnotPlot by T. Pogmore, for both inversion experiments and knot relaxations.

\section{Notation and basic properties of sphere inversions}

\label{Basics}

For a given center $\mathbf{p}$ and radius $r$, let $\rho[\p,r]$
denote inversion in the sphere of radius $r$ centered at $\p$. Suppose $K$
is a polygon, and $\p$ is not one of the vertices of $K$.  Then $\rho[\p,r](K)$, the image of $K$ under the inversion
map, is a knot in $\R^3$ equivalent to $K^*$. Let $\hat\rho[\p,r](K)$ denote the polygon obtained by inverting the
vertices of $K$ and reconnecting them with line segments. We sometimes
suppress the $[\p,r]$ and just write $\rho(K)$ and $\hat\rho(K)$, more
briefly $\widehat{K}$.

The map $\rho[\p,r]: \R^{3}-\{\p\}\to
\R^{3}-\{\p\}$ is a conformal homeomorphism. (The derivative
$D_{\x}\rho$ is a similarity.)

If $\p$ and $r$ become infinite at the same rate, then the polygons $\rhat(K)$ approach $K^*$.  For example, let $\p = (r,0,0)\in \R^3$.  Then for any $(x,y,z)\in \R^3$,
\[
\lim_{r\to\infty} \rho(x,y,z) = (-x,y,z).
\]

The composition of two inversions (with different radii) with the same center $\p$ is a dilation from $\p$.  When $\p \neq \0$, the map is an affine isomorphism, the composition of a translation and a dilation from $\0$.

Any finite composition of sphere inversions is called a {\em M\"obius transfmormation}.  If a M\"obius transformation $\mu : \R^3 \cup \infty \to \R^3 \cup \infty$ happens to take $\infty \to \infty$ (i.e. takes $\R^3 \to \R^3$), then $\mu | \R^3$ is an affine isomorphism, in particular a composition of a translation, rotation, similarity, and perhaps reflection in a plane.

\section{Knot types are determined by sphere complements}

\label{ReduceToSpheres}

Suppose $K$ is a given polygonal knot. We want to bound the number of distinct
knot types that can arise as $\hat\rho(K)$, considering all possible spheres
of inversion. The first observation is that the knot type of $\hat\rho(K)$
depends only on the choice of center $\mathbf{p}$, not the radius of the
sphere. As noted in the previous section, two inversions in the same center
with different radii differ by a Euclidean similarity, which preserves knot
type. We next want to show that the possible centers of inversion fall into a
relatively small number of equivalence classes.

Fix the radius of inversion to be $1$, and suppose we have two centers of
inversion, $\mathbf{p}$ and $\mathbf{q}$. Let $\al(t),\;t=0\ldots1$, be a
path (missing $K$) with $\al(0)=\mathbf{p}$ and $\al(1)=\mathbf{q}$.
This gives a homotopy of inversion maps and a homotopy of polygons $\hat
\rho[\al(t),1](K)$. If none of the polygons $\hat\rho[\al(t),1](K)$ is
self-intersecting, in the sense that two non-adjacent edges intersect, then
the knot types of all $\hat\rho[\al(t),1](K)$, in particular $\hat
\rho[\mathbf{p},1](K)$ and $\hat\rho[\mathbf{q},1](K)$, are the same. Of
course, sometimes $\hat\rho[\al(t),1](K)$ is self-intersecting, which is
why the process being studied can produce many knot types.

When can $\hat\rho(K)$ be singular?. Suppose $E$ and $F$ are two nonadjacent
edges of $\hat\rho(K)$ that meet. Then $E$ and $F$ must be co-planar, in
particular the four vertices $\rho(\mathbf{w}), \rho(\mathbf{x}),
\rho(\mathbf{y}), \rho(\mathbf{z})$ of $E$ and $F$ must be co-planar.

If the four points $\rho(\mathbf{w}), \rho(\mathbf{x}),
\rho(\mathbf{y}), \rho(\mathbf{z})$ are contained in a unique plane, then
the four points $\mathbf{w},\mathbf{x},\mathbf{y},\mathbf{z}$ lie on a unique
round 2-sphere [or plane] $S$; and the center of inversion, $\mathbf{p}$, must also lie on S. (i.e. coplanar with $\mathbf{\infty} \rightleftarrows$ cospherical with $\mathbf{p}$.) So long as the path $\al(t)$ from
$\mathbf{p}$ to $\mathbf{q}$ can be chosen to not intersect $S$,
we can conclude that $\hat\rho[\mathbf{p},1](K)$ and $\hat\rho[\mathbf{q}%
,1](K)$ are the same knot type.

The other possibility for a 4-tuple of vertices is that the points
$\rho(\mathbf{w})$, $ \rho(\mathbf{x})$,  $\rho(\mathbf{y})$,
 and  $\rho(\mathbf{z}) $ do not determine a unique plane, in which case they are
colinear. Then $\mathbf{w},\mathbf{x},\mathbf{y},\mathbf{z}$ lie on a some
circle. Any path $\al(t)$ that misses the 2-spheres described above can be
wiggled slightly to miss any of these circles as well.

If $\mathbf{p}$ lies on one of the 2-spheres and $\hat\rho[\mathbf{p},1](K)$
is nonsingular (so we \underbar{do} want to reckon with its knot type) then $\hat
\rho[\mathbf{p},1](K)$ has the same knot type as for centers near $\mathbf{p}$
on either side of the sphere.

We conclude the following:

\begin{theorem}  \label{BoundSpheres}
The number of knot types that can arise as $\hat
\rho[\mathbf{p},r](K)$, over all $r \in\mathbb{R}$ and $\mathbf{p}
\in\mathbb{R}^{3} - \{\text{vertices of K}\}$ is bounded by the number of complementary domains of the
system of round 2-spheres [and planes] consisting of the unique spheres [or planes] determined by
4-tuples of vertices of non-adjacent edges of $K$. There are at most $ n(n-3)/2$
such spheres, where $n$ is the number of vertices of $K$.
\end{theorem}

\section{Counting complementary domains of spheres}

\label{Domains}

Our goal in this section is to count the number of complementary regions for a
collection of $m$  round two-spheres (which might intersect in various ways) in $\mathbb{R}^{3}$. \ We first give an elementary proof of an upper bound, and then establish a formula (Theorem \ref{exactformula})  for the exact number, using work of
Ziegler and \v{Z}ivaljevi\'{c} \cite{ZZ}. \ From Theorem \ref{exactformula}, we shall see that the bound given in Theorem \ref{maximal} is sharp in the sense that it is attained by any collection of spheres that intersect generically (see Section \ref{genericspheres}).

\begin{definition}
A \emph{round sphere }$S$ in $\mathbb{R}^{3}$ is any set of the form
\begin{equation*}
S=\{(x,y,z)\mid (x-x_{0})^{2}+(y-y_{0})^{2}+(z-z_{0})^{2}=r_{0}^{2}\},\text{
}r_{0}>0
\end{equation*}%
A \emph{round circle }is any non-empty transverse intersection of a round
sphere and a plane.
\end{definition}
If $S_1$ and $S_2$ are round spheres, then $S_1 \cap S_2$ must be either empty, one point, one round circle, or  $S_1 = S_2$.

\begin{proposition} \label{CirclesInSphere}
Let $\{C_{1},\ldots ,C_{k}\}$ be a collection of $k$ round circles in a
round sphere $S$. \ Then the complement $S\setminus \cup C_{i}$ has at most $%
k^{2}-k+2$ components.
\end{proposition}

\begin{proof}
The result clearly holds for $k=1$. \ Assume it holds for $k-1$ circles. \ The
circle $C_{k}$ intersects $\cup_{i=1}^{k-1}C_{i}$ in at most $2(k-1)$ points,
which separate $C_{k}$ into at most $2(k-1)$ arcs. \ Thus, since each arc of
$C_{k}$ can separate at most one complementary region of $S\setminus\cup
_{i=1}^{k-1}C_{i}$ in two, the number of complementary regions of
$S\setminus\cup_{i=1}^{k}C_{i}$ is at most $\left[  (k-1)^{2}-(k-1)+2\right]
+2(k-1)=k^{2}-k+2$.
\end{proof}

\begin{theorem}\label{maximal}
Let $\{S_{1},\ldots,S_{m}\}$ be any collection of $\ m$ round two-spheres in
$\mathbb{R}^{3}$. \ Then $\mathbb{R}^{3}\setminus\cup_{i=1}^{m}S_{i}$ has at
most
\[
m^{3}/3-m^{2}+8m/3=2\binom{m}{3}+2m
\]
components.
\end{theorem}

\begin{proof}

The theorem is clearly true for $m=1$. \ Assume it holds for $m-1$ two-spheres.
\ By the previous result the two-sphere $S_{m}$ is cut into at most
$(m-1)^{2}-(m-1)+2$ regions by the round circles $S_{i}\cap S_{m}$, for
$i=1,\ldots,m-1$. \ (Here some $S_{i}\cap S_{m}$ may be points or empty.)
\ Thus $\mathbb{R}^{3}\setminus\cup_{i=1}^{m}S_{i}$ has at most
\[%
\begin{array}
[c]{c}%
\frac{(m-1)^{3}}{3}-(m-1)^{2}+\frac{8}{3}(m-1)+(m-1)^{2}-(m-1)+2\\
=\frac{m^{3}}{3}-m^{2}+\frac{8}{3}m=2\binom{m}{3}+2m
\end{array}
\]
components.
\end{proof}

We shall use  \cite[Corollary 2.8 and
Theorem 2.7]{ZZ} to obtain a general formula for the cohomology groups of 
$\mathbb{R}^{3}\setminus\cup_{i=1}^{m}S_{i}$.  \ Once again, let $\{S_{1},\ldots,S_{m}\}$ be any collection of $\ m$ round two-spheres in
$\mathbb{R}^{3}$. \ It is clear that the intersection of any collection of some of these spheres is either some
$S_{i}$, a round circle, two points, a single point, or empty. \ Let $P$ be
the partially ordered set (poset) of connected components of these
intersections, ordered by reverse inclusion. \ The \emph{order complex}
$\Delta(P)$ is the simplicial complex with vertices the elements of this
partially ordered set, and simplices corresponding to ordered linear chains.
\ We include the empty set as maximal element, denoted $\widehat{1}$.

\begin{example}
Suppose we have four two-spheres $S_{1},\ldots,S_{4}$ with $S_{1}\cap
S_{2}$, $S_{1}\cap S_{3}$, and $S_{2}\cap S_{3}$ circles.  Denote the circles $C_{12}$,  $C_{13}$, and $C_{23}$ respectively.  Suppose further that $S_{1}\cap S_{2}\cap S_{3}$
  is exactly two points,  $S_{1}\cap S_{4}=S_{2}\cap S_{4}=\emptyset\,$, and $S_{3}\cap
S_{4}$ is a single point.  Then the poset has twelve elements, and the associated
order complex has twelve vertices and is  of dimension three.
\end{example}

Let $P_{0}$ be the subposet of all elements which are spheres, including the
maximal element $\widehat{1}$.

Then, letting $\simeq$ denote homotopy equivalence and $P_{<p}$ denote the
subposet of elements less than $p$, Theorem 2.7 of \cite{ZZ} gives
\[
\cup S_{i}\simeq\Delta(P_{<\widehat{1}})\vee(\Delta(P_{<p})\ast S^{d(p)}%
\]
where the wedge product $\vee$ is taken over all $p\in P_{0}\setminus
\widehat{1}$, $\ast$ denotes join, and $S^{d(p)}$ is a sphere of dimension
equal to the geometric dimension $d\left(  p\right)  $ of $p$. \ Now by
Alexander duality, Corollary 2.8 of \cite{ZZ} gives

\begin{theorem} \label{exactformula}
The cohomology groups of the complement of a union of round two-spheres are
given by%
\[
\widetilde{H}^{i}(\mathbb{R}^{3}\setminus\cup S_{i})\cong\bigoplus_{p\in
P,d(p)\neq0}\widetilde{H}_{1-d(p)-i}(\Delta(P_{<p}))\;.
\]

\end{theorem}

Recall that the number of path components of a space is equal to the rank of the zero-th unreduced cohomology group of the space.

In the above {\bf Example}, $\{p\in P\mid d(p)\neq0\}=\{S_{1},S_{2},S_{3},S_{4}%
,C_{12},C_{13},C_{23},\widehat{1}\}$ and terms contributing to homology in the
direct sum are (i) $\widetilde{H}_{0}(\Delta(P_{<C_{ij}}))$ where
$\Delta(P_{<C_{ij}})$ is two points for each $C_{ij}$, (ii) $\widetilde{H}%
_{2}(\Delta(P_{<\widehat{1}}))$, and (iii) one copy of the integers for each
$S_{i}$ (by convention $\widetilde{H}_{-1}(\emptyset)\cong\mathbb{Z}$). 
 So
contributions to the rank of $\widetilde{H}^{0}(\mathbb{R}^{3}\setminus\cup
S_{i})$ are one from each $C_{i}$, one for each $S_{i}$, and one for
$\widetilde{H}_{2}(\Delta(P_{<\widehat{1}}))$. \ To see the last, note that in
$\Delta(P_{<\widehat{1}})$ there is a hexagon formed by six edges (${C}_{12}$ to $S_{1}$ to $C_{13}$ to $\ldots$to $S_{3}$ to $C_{13}$). \ This
hexagon is suspended at two points $P_{1}$ and $P_{2}$. \ Finally, there are
two additional vertices $Q$ and $S_{4}$, and two edges $S_{3}Q$ and $QS_{4}$.
\ Thus in the example the complement of the four round two-spheres has nine
path components.

\subsection{Collections of spheres that intersect generically} \label{genericspheres}

Let us compare the result just proved with our earlier consideration
of the maximal case. When the $m$ round two-spheres intersect generically,
one has
\begin{itemize}

\item The original m copies of $S^{2}$.
\item As double intersections, $\binom{m}{2}$ copies of round circles.
\item As triple intersections, $\binom{m}{3}$ copies of $S^{0}$ (i.e. $2\binom{m}{3}$ points).
\end{itemize}

In order to determine the number of path components, we need to understand
$\Delta(P_{<\widehat{1}})$. \ A direct count of simplices shows that the euler
characteristic is
\[
2\binom{m}{3}-\binom{m}{2}+m
\]
Clearly $\Delta(P_{<\widehat{1}})$ is simply connected (it has the homotopy
type of a wedge of two-spheres), so the second betti number $b_{2}
(\Delta(P_{<\widehat{1}}))=2\binom{m}{3}-\binom{m}{2}+m-1$. \ Thus the number
of path components of $\mathbb{R}^{3}\setminus\cup S_{i}$ in the generic case
is
\[
m+\binom{m}{2}+2\binom{m}{3}-\binom{m}{2}+m-1+1=2\binom{m}{3}+2m
\]
which checks with our calculation of the maximum number earlier.  Thus the generic case is the maximal case.

\subsection{Apply to spheres coming from  polygons} 
\label{spheresfrompolygons}


\begin{corollary} \label{BoundKnots}
The number of knot types obtainable by inverting a given n-edge knot is at most
\[
\frac{1}{24}\left(n^6-9n^5+21n^4+9n^3-22n^2-96n\right)
\]
\end{corollary}

\begin{proof}
By Theorem \ref{BoundSpheres}, the number of knots is bounded by the number of complementary domains of (at most) $\frac{n(n-3)}{2}$ spheres.  Use this value for $m$ in the bound from Theorem \ref{maximal}.
\end{proof}

\section{Exponentially many $n$-segment knot types}

\label{ExponentiallyMany}

\begin{theorem} \label{NumberOfKnotTypes}
The number of knot types representable by polygons with $n$ or fewer edges is exponential in $n$, in particular $\geq \frac{1}{12}(\sqrt{2}^n - 4)$.

\end{theorem}

\begin{proof}
Ernst and Sumners showed \cite{ErnstSumners} that the number of prime knot
types (distinguishing mirror images) that can be represented with diagrams
having $\leq q$ crossings is at least $(1/3)(2^{q-2}-1)$, so the total of all
$q$-crossing knot types is even greater. On the other hand, Negami showed \cite{Negami} that if a given
knot type $[K]$ can be realized with some $q$-crossing knot, then $[K]$ can be
realized as a polygon with $\leq2q$ edges. Thus the number of $n$-edge knot
types is at least as large as the number of $(n/2)$-crossing knots types, in particular exponential in $n$. Substituting $q=n/2$ into the Ernst-Sumners bound gives the lower bound $\frac{1}{12}(\sqrt{2}^n - 4)$
\end{proof}

\begin{remark}
Strictly speaking, we should distinguish the cases where $n$ is odd vs. even in the above calculation, since $q=n/2$ must be an integer number of crossings.  So for odd $n$, we should substitute $q=(n/2)-1/2$ in the Ernst-Sumners formula.  However, this bound only talks about prime knots.  If we include composite knots in giving a lower bound for the total number of $q$ crossing knot types, then by the time we get to q=20 or more crossings,   the total number of knot types is indeed larger than $(1/3)(2^{q-2}-1)$.
\end{remark}

\section{The whole group of M\"obius transformations} \label{MobiusGroup}

  Suppose $\mu$ is any M\"obius transformation such that none of the vertices of $K$ is taken by $\mu$ to $\infty$.  Then we can define $\hat\mu(K) \subset \R^3$, as before, by connecting the points $\mu(v_1), \ldots, \mu(v_n)$ with lines.

\begin{theorem} \label{BoundMobius}
The number of knot types  that can arise as $\hat\mu(K)$ is at most twice the number that can arise as $\rhat (K)$ for single inversions.  Specifically, if a knot is obtainable as $\hat\rho(K)$ then its mirror image is obtainable by an additional inversion; but otherwise, composing inversions does not discover additional knot types.\end{theorem}

\begin{proof}
Let $\p = \mu^{-1}(\infty)$.  Note $\p$ is not one of the vertices of $K$.

If $\p=\infty$, then $\mu|\R^3$ is an affine isomorphism, so $\hat\mu(K)$ = $\mu(K)$, which is the same knot type as $K$ or $K^*$.

If $\p \neq \infty$, then let $\rho$ be inversion in a sphere of some radius centered at $\p$.  The map $\mu \circ \rho$ is a M\"obius transformation sending $\infty \to \infty$, hence an affine isomorphism of $\R^3$.  Thus $\hat\mu(K) = (\mu \circ \rho)(\rhat(K))$ is of the same knot type as $\rhat(K)$.  So the  knot types that can arise as  $\hat\mu(K)$ are precisely those that  arise as some $\rhat(K)$ and their mirror images.

\end{proof}

\section{Conclusion}

The number of knot types obtainable by inverting a given n-edge knot is bounded (Corollary \ref{BoundKnots}) as
\[
\text{Number obtainable } \leq \frac{1}{24}\left(n^6-9n^5+21n^4+9n^3-22n^2-96n\right)\;
\]

If we allow the whole M\"obius group, i.e. allow finite compositions of inversions, then, 
from Theorem \ref{BoundMobius},
the upper bound  is at most doubled.

On the other hand, from Theorem \ref{NumberOfKnotTypes}, the total number of $n$-edge knot types is bounded below as
\[
\text{All n-edge knots } \geq \frac{1}{12}(\sqrt{2}^n - 4)\;.  
\]

Thus for $n=72$ edges and beyond, no given $n$-edge polygon can generate all $n$-edge knot types by single inversion.  And for $n \geq 75$, no polygon can be universal under the whole M\"obius group.

\clearpage

\section{Figures}

\clearpage

\vspace{0.5in}

\clearpage

\vspace{0.5in}

\clearpage


\clearpage

\section{Acknowledgements}
Please see the comments in Section \ref{ackn}.

\end{document}